\documentclass{article}

\usepackage{amsmath, amsthm, amssymb}
\usepackage{graphicx}
\usepackage{verbatim}
\usepackage{natbib}
\usepackage{caption}
\usepackage{subcaption}
\usepackage{fancyvrb}
\usepackage{enumerate}
\usepackage{relsize}
\usepackage{mathtools}
\usepackage{enumitem}

\usepackage{hyperref}
\usepackage[margin=1.5in]{geometry}
\hypersetup{colorlinks,citecolor=blue,urlcolor=blue,linkcolor=blue}

\usepackage{stefan_tex}
\graphicspath{{./figures/}}

\DeclarePairedDelimiter\abs{\lvert}{\rvert}
\DeclarePairedDelimiter\norm{\lVert}{\rVert}
\DeclarePairedDelimiter\set{\{}{\}}
\theoremstyle{plain}
\newtheorem{prop}{Proposition}

\newtheorem{theo}[prop]{Theorem}

\theoremstyle{definition}

\theoremstyle{remark}

\author{David A.~Hirshberg  \and Stefan Wager}
\date{Stanford University}
\title{Debiased Inference of Average Partial Effects \\ in Single-Index Models}

\begin{document}

\maketitle

\begin{abstract}
We propose a method for average partial effect estimation in high-dimensional single-index models that is $\sqrt{n}$-consistent and asymptotically unbiased given sparsity assumptions on the underlying regression model. This note was prepared as a comment on \citet{wooldridge2018inference}, forthcoming in the Journal of Business and Economic Statistics.
\end{abstract}

\paragraph{Introduction}
There has recently been a considerable amount of interest in developing methods for statistical
inference in high-dimensional regimes with more covariates
than data points \citep{athey2016approximate,belloni2013program,javanmard2014confidence,
van2014asymptotically,zhang2014confidence}.
\citet{wooldridge2018inference} build on this literature, and propose a new method for inference
about average partial effects in high-dimensional probit models; they then extend their
approach to non-linear panels with correlated random effects \citep{wooldridge2010econometric}.
This is a valuable result, with many potential application areas.

In order to achieve $\sqrt{n}$-consistent inference, however, the method studied by \citet{wooldridge2018inference}
requires a ``soft'' beta-min condition that asymptotically
rules out regularization bias from model selection. And, as argued by  \citet{belloni2014inference}, this
type of approach may be vulnerable to confounding when there are features that are highly correlated
with the focal variable and have weak but non-zero effects on the outcome.

In this comment, we discuss an alternative approach to average partial effect estimation that avoids using a
beta-min-style assumption by explicitly accounting for the correlation structure of the features. Qualitatively,
our approach is related to both the double-selection principle \citep{belloni2014inference,chernozhukov2016double}
and the idea of modeling or balancing the propensity score for average treatment effect estimation
\citep{athey2016approximate,farrell2015robust,robins1}. Formally, we apply the debiasing idea of
\citet{javanmard2014confidence} to a linearization of the original problem. We prove that our approach
allows for $\sqrt{n}$-consistent inference under assumptions that are more closely in line with the broader literature, and
show in simulations that this approach can be more robust that of \citet{wooldridge2018inference} in the presence of confounding.

We focus on inference in high-dimensional single-index models.
We observe $n$ independent and identically distributed samples
\smash{$(X_i, \, Y_i) \in \RR^p \times \cb{0, \, 1}$}, where $p$ may be much larger than $n$,
and seek to estimate the average partial effect (APE) $\tau_j$  for
some $j \in 1, \, ..., \, p$,
\begin{equation}
\label{eq:APE}
\EE{Y\cond X = x} = \Psi\p{x^T \theta}, \ \ \tau_j = \EE{\frac{d}{d X_j} \EE{Y \cond X}} = \theta_j \EE{\psi\p{X^T \theta}},
\end{equation}
where $\Psi(\cdot)$ is a link function with derivative $\psi(\cdot)$.
The simplest specification considered by \citet{wooldridge2018inference} corresponds to \eqref{eq:APE} with a probit
link, i.e., $\Psi(\cdot) = \Phi(\cdot)$ for the standard Gaussian cumulative distribution function $\Phi(\cdot)$.
In this note, we do not consider the richer class of panel models discussed in
\citet{wooldridge2018inference}, and simply focus on the i.i.d. case.

\paragraph{Background: Single-Index Models and Non-Linear Estimation}

The main difficulty of this problem relative to existing results on debiased estimation is that $\tau_j$ as specified in \eqref{eq:APE} in a non-linear
functional of $\theta$. In contrast, the problem of estimating average effect of a binary treatment with high-dimensional confounding
is a linear problem, and so can be approached with a more standard debiasing approach \citep{athey2016approximate}.
Other papers that discuss the problem of non-linearities in high-dimensional inference include
\citet{van2014asymptotically} and \citet{belloni2013program}; see \citet{wooldridge2018inference} for a discussion.

It is interesting to consider why the task of APE estimation, as framed here, results in a non-linear problem. After
all, when written in terms of the conditional response surface $m(x) = \EE{Y \cond X = x}$, the APE is linear in $m$. For
example, if we write down a conditionally linear model for $m$ as below,
then we can write $\tau$ as a simple weighted average of $Y$:
\begin{equation}
\label{eq:cond_lin}
\text{if } \, m(x) = \mu\p{x_{-j}} + x_j \delta\p{x_{-j}}, \, \text{ then } \, \tau = \EE{\frac{\p{X_j - \EE{X_j \cond X_{-j}}} Y}{\Var{X_j \cond X_{-j}}}},
\end{equation}
where $x_{-j}$ denotes the $p-1$ dimensional vector obtained by removing the $j$-th entry of $x$. More generally,
under simple regularity conditions, we can express $\tau$ in terms of the density of
$X_j$ conditionally on $X_{-j} = x_{-j}$, denoted $f_j(\cdot; x_{-j})$ \citep{powell1989semiparametric}:
\begin{equation}
\label{eq:parts}
\tau = \EE{-\frac{d}{d z} \cb{\log\p{f_j\p{z; \, X_{-j}}}}_{z = X_j} \, Y}.
\end{equation}
In both cases, the representations in \eqref{eq:cond_lin} and \eqref{eq:parts} can guide semiparametrically efficient estimation
of $\tau$, either via debiased estimation as considered here \citep{hirshberg2017balancing} or via plug-in estimation using
an appropriately chosen orthogonal moments construction \citep{chernozhukov2016locally,chernozhukov2018double}.

However, a key aspect of both the conditionally linear model \eqref{eq:cond_lin}
and the fully generic model underlying \eqref{eq:parts}
is that the underlying model class for $m(\cdot)$ is convex, and convexity plays a key role in enabling practical inference
about linear functionals \citep{armstrong2015optimal,donoho1994statistical,hirshberg2017balancing}.
Here, conversely, the class of functions
$m_\theta(x) = \Psi(x \theta)$ is not convex in $m$-space, and so the machinery used to prove semiparametric efficiency
in  \citet{chernozhukov2016locally,chernozhukov2018double} or \citet{hirshberg2017balancing} is not immediately available. 

\paragraph{Debiasing in Single-Index Models}

We now return to our main focus, that is debiased inference of average partial effects as defined in \eqref{eq:APE}.
As in \citet{chernozhukov2016locally,chernozhukov2018double} and \citet{hirshberg2017balancing}, we study estimators that start with a
parameter estimate \smash{$\htheta$}, and then debias the plug-in estimator for $\tau_j$ based on \smash{$\htheta$} using a weighted
average of residuals:\footnote{Our proposed estimator will
also use cross-fitting \citep{chernozhukov2016double}, but we suppress this notation here for conciseness.}
\begin{equation}
\label{eq:debias}
\htau_j = \frac{1}{n} \sum_{i = 1}^n \p{\htheta_j \psi\p{X_i^T \htheta} + \hgamma_i \p{Y_i - \Psi\p{X_i^T \htheta}}},
\end{equation}
In contrast, the approach of \citet{wooldridge2018inference} takes on a markedly different functional form.
Their method first gets an estimate $\htheta$ via $L_1$-penalized quasi-maximum likelihood estimation,
and selects a signal set \smash{$A \subset \cb{1, \, ..., \, p}$} that contains the non-zero entries of \smash{$\htheta$} along with some
pre-determined variables of interest. They then obtain a corrected estimator \smash{$\ttheta$}, where
\smash{$\ttheta_j$} for $j \in A$ is obtained via a generalization of the debiased lasso of \citet{javanmard2014confidence},
while \smash{$\ttheta_j = 0$} for \smash{$j \not\in A$}.
Finally, they conclude with a plug-in step, 
\begin{equation}
\label{eq:debiasWZ}
\htau_j^{WZ} = \frac{1}{n} \sum_{i = 1}^n \ttheta_j \psi\p{X_i^T \ttheta}.
\end{equation}
One key difference relative to existing results
is that the method only debiases the set of covariates $A$ that are useful for predicting $Y$ (along with
a deterministic pre-determined set), rather than correcting the sample
Hessian over all covariates as in \citet{athey2016approximate} or \citet{javanmard2014confidence}.
This may make the procedure less robust in cases where some
variables with weak but non-zero signals are strongly predictive of $X_j$ \citep{belloni2014inference},
and is reflected in the beta-min-style condition discussed above.

\citet{hirshberg2017balancing} showed that estimators like \eqref{eq:debias} for linear functionals of $m(x)=\EE{Y_i \mid X_i=x}$
are semiparametrically efficient with considerable generality when the weights \smash{$\hgamma_i$} solve a minimax problem: 
minimizing the maximal conditional-on-$X$ MSE of our estimator \eqref{eq:debias} over a set of plausible regression error functions $\hm - m$.
This result builds on key contributions of \citet{donoho1994statistical} and \citet{chernozhukov2016locally}.
Although deriving these minimax weights \smash{$\hgamma_i$} for our present problem would be
difficult because our estimand is nonlinear in $\theta$, we can still find weights that solve a first-order approximation to this minimax problem. 
In particular, if we have an estimator $\smash{\htheta}$ that we believe to be accurate in $\ell_1$ norm, we choose the weights
\begin{equation}
\label{eq:balance}
\begin{split}
&\hgamma_i = \argmin_\gamma\cb{I^2\p{\gamma; \, \htheta} + \frac{1}{n^2} \Norm{\gamma}_2^2}, \\ 
&I\p{\gamma; \, \htheta} = \Norm{\frac{1}{n} \sum_{i = 1}^n  \p{\htheta_j \psi'\p{X_i^T \htheta} X_i + \psi\p{X_i^T \htheta} e_j  - \gamma_i \psi\p{X_i^T \htheta} X_i} }_{\infty}.
\end{split}
\end{equation}
We derive this problem by Taylor expansion of our estimator's error. 
Writing the sample-average version of our estimand as \smash{$\htau_j^* = n^{-1} \sum_{i = 1}^n \theta_j \psi\p{X_i \theta}$},
 we can characterize the error of the estimator \eqref{eq:debias} as follows, subject to a few conditions stated in the theorem below.
\begin{equation}
\label{eq:error-expansion}
 \begin{aligned}
 \htau^*_j - \htau_j
 &= \frac{1}{n} \sum_{i = 1}^n \p{\theta_j \psi\p{X_i^T \theta} - \htheta_j \psi\p{X_i^T \htheta} - \hgamma_i \p{Y_i - \Psi\p{X_i^T \htheta}}} \\
 &= \p{\frac{1}{n} \sum_{i = 1}^n \htheta_j \psi'\p{X_i^T \htheta} X_i^T + \psi\p{X_i^T \htheta} e_j^T -  \hgamma_i \psi\p{X_i^T \htheta} X_i} \p{\theta - \htheta} \\
 &\ \ \ \ \ \ \ \ - \frac{1}{n} \sum_{i = 1}^n \hgamma_i \p{Y_i - \Psi\p{X_i^T \theta}} + rem_n, \ rem_n \in o_P\p{\lVert \htheta - \theta \rVert_2^2}. \\
 \end{aligned}
\end{equation}
Our optimization problem \eqref{eq:debias} is chosen to make the first two terms above small. Its first term is the square of the H\"older's inequality bound on our first term above for $\lVert\theta - \htheta\rVert_1 = 1$ and its second is the mean square of the second term above when $\Var{Y_i \mid X_i}=1$. 
The result below establishes formal conditions under our proposed estimator provides consistent $1/\sqrt{n}$-scale estimates and asymptotically valid 
confidence intervals for $\tau_j$. To simplify the statement and proof of this result, we make the impossible assumption that
we have a deterministic-yet-consistent pilot estimator $\htheta$ of $\theta$. Analogous results can be proven 
for a pilot estimator $\htheta$ independent of $(X_i, Y_i)_{i \le n}$ defined on an auxilliary sample
and which satisfies with probability tending to one the same properties we require of our deterministic estimator
sequence $\htheta$; moreover, we can use cross-fitting to avoid efficiency loss from sample splitting
\citep{chernozhukov2016double}.

\begin{theo}
 \label{theo:main}
	Suppose that we observe $(X_{i,n},Y_{i,n})_{i \le n}$ iid with $Y_{i,n} \in \mathbb{R}$ and $X_{i,n} \in [-1,1]^{p_n}$ for $\log(p_n) = o(n)$
	and that $\EE[]{Y_{i,n} \mid X_{i,n}=x}=\Psi(x^T \theta_n)$ for some link function $\Psi$ with 3 bounded derivatives and $\theta_{j,n} = O(1)$.
    Suppose, in addition, that we have a deterministic estimator sequence $\htheta_n$ that satisfies $\norm{\htheta_n - \theta_n}_{1} = o(1)$
	and $(\htheta_n - \theta_n)^T A (\htheta_n - \theta_n) = o(n^{-1/2})$ for $A \in \mathcal{A}_n = \set{  E[ X_{i,n} X_{i,n}^T],\ E[ X_{i,n} e_j^T + e_j X_{i,n}^T] }$.
	Then in terms of $\psi = \Psi'$, define $\gamma^{\star}_{n}(x) = \psi(x^T \htheta_n) x^T g(\htheta_n)$ for
	\begin{equation}
	\label{eq:riesz-weights}
	g(\theta) = \EE[]{\psi(X_i^T \theta)^2 X_i X_i^T}^{-1} \EE[]{\theta_{j} \psi'(X_i^T \theta) X_i + \psi(X_i^T \theta)e_j}
	\end{equation}
    If $\norm{g(\htheta_n)}_{1} = o(n^{1/2})$ and $\norm{\gamma^{\star}_{n}}_{\infty} = O(1)$
	and we consider $\htau_j$ as in \eqref{eq:debias} with $\htheta_n$ as above and $\hgamma$ as in \eqref{eq:balance},
	we have the asymptotic characterization
	\begin{equation}
	\label{eq:asymptotic-linear}
	\begin{split}
	&\htau_j - \tau_j = n^{-1}\sum_{i=1}^n \iota_n(X_i,Y_i) + o_p(n^{-1/2}) \ \text{ for }\\
	&\iota_n(x,y) =  \theta_{j,n}\psi(x^T \theta_n) - \tau_j  +  \gamma_n^{\star}(x)\p{y - \Psi(x^T \theta_n)}.
	\end{split}
	\end{equation}
\end{theo}

This asymptotic characterization implies that for $V_n = E \iota_n(X_i,Y_i)^2$, $\sqrt{n}(\htau_j - \tau_j)/V_n^{1/2}$ will be asymptotically normal with variance one,
justifying inference as usual.  We make essentially three assumptions. Our first assumption, just
as in \citet{wooldridge2018inference}, is the correctness of a parametric model $E[Y \mid X=x]=\Psi(x^T \theta_n)$ with $\theta_{j,n}=O(1)$;
our second one is that our pilot estimator $\htheta_n$ is $\ell_1$-consistent and
satisfies additional consistency properties discussed below and in our proof;
and our third, a type of identifiability condition, requires that we observe adequate variation 
in $X_i$ in the two directions in which differences between $\theta$ and our pilot estimate $\htheta$ 
will, if uncorrected, result in significant bias in our estimate of $\tau_j$: the direction $e_j$ and the direction $E[\psi'(X_i^T \htheta) X_i]$.

We end our discussion with some brief comments about the assumptions used to prove
Theorem \ref{theo:main}. First, as emphasized earlier, our proof does not rely on recovering
the support of $\theta$, and thus we do not require any form of beta-min condition (i.e., the non-zero
entries of $\theta$ are allowed to be very close to 0). This may make our result more robust in
the presence of weak signals.

We make several high-level assumptions about the behavior of the pilot estimator
\smash{$\htheta$}. If we are willing to assume that the maximal eigenvalues of both matrices in
the set $\mathcal{A}_n$ are bounded uniformly in $n$, then it is sufficient to assume that
\smash{$\norm{\htheta-\theta}_{1} = o(1)$} and \smash{$\norm{\htheta - \theta}_{2} = o(n^{-1/4})$}.
It is well known that if $\theta$ is $k$-sparse for some \smash{$k \ll \sqrt{n}/\log(p)$},
then we can obtain estimators \smash{$\htheta$} that satisfy these bounds with high probability
using different variants of $\ell_1$-penalized regression \citep{hastie2015statistical}.
The implicit sparsity assumption \smash{$k \ll \sqrt{n}/\log(p)$} is substantially weaker than the corresponding
assumption made in Theorem 4.1 of \citet{wooldridge2018inference}, namely \smash{$k \ll (\sqrt{n}/\log(p))^{2/3}$}.

\paragraph{Numerical Experiment}

An in-depth empirical investigation is beyond the scope of this comment. However, we illustrate the behavior
of our new estimator of average partial effects via augmented minimax linear estimation in a simple simulation
experiment. Just as in \citet{wooldridge2018inference}, we study the binary outcome case; however, we use a logistic rather
than probit link function because penalized logistic regression is readily accessible via the \texttt{R}-package
\texttt{glmnet} \citep{friedman2010regularization}. We also compare our method to the natural logistic variant of
the method proposed by \citet{wooldridge2018inference}.\footnote{Replication files for all experiments are available at
\url{github.com/swager/amlinear}, in the folder \texttt{debiased\_single\_index\_experiments}. For convex
optimization, we use the \texttt{R} package \texttt{CVXR} \citep{CVXR}.}

In all our experiments, we generate data as follows, where $x \in \RR^p$:
\begin{equation}
\PP{Y \cond X = x} = \frac{1}{1 + e^{-x \cdot \theta}}, \ \ \theta_1 = -\frac{1}{10}, \ \ \theta_j = \frac{20}{\p{5 + j}^2} \ \eqfor  \ j = 2, \, \ldots, \, p,
\end{equation}
and seek to estimate $\tau_1$, the average partial effect with respect to the first feature.
We consider both a setting with uncorrelated features, \smash{$X \sim \nn\p{0, \, \ii_{p \times p}}$},
and correlated features:
\begin{equation}
X_{2:p} \sim \nn\p{0, \, \ii_{(p - 1) \times (p - 1)}}, \ \ X_1 \cond X_{2:p} \sim \nn\p{\sqrt{\frac{1}{30}} \sum_{j = 11}^{20} X_j, \, \frac{2}{3}}.
\end{equation}
We varied the sample size $n$, and always set $p = 2n$.

\begin{figure}[t]
\begin{center}
\begin{tabular}{cc}
\includegraphics[width=0.45\textwidth]{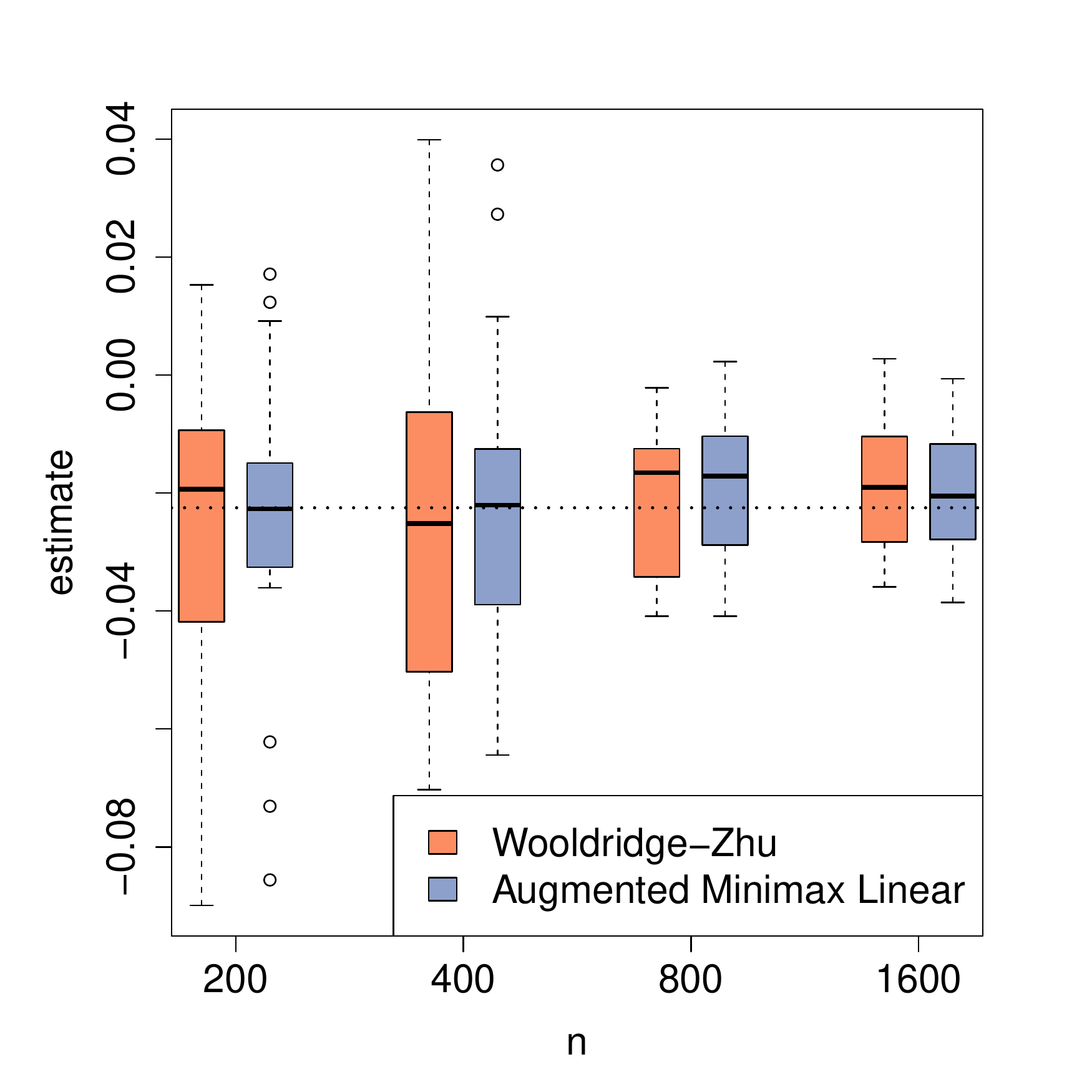} &
\includegraphics[width=0.45\textwidth]{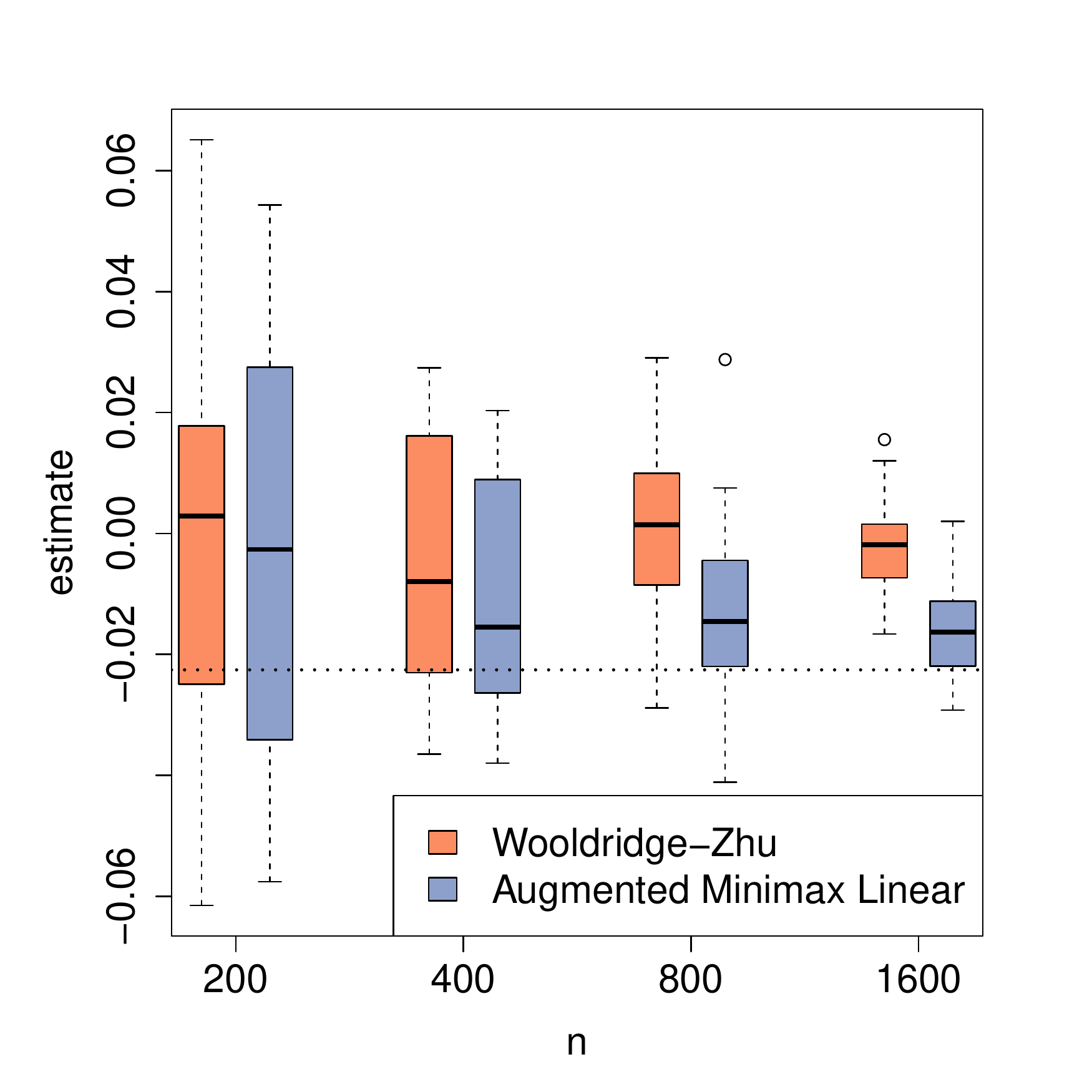} \\
uncorrelated $X$ & correlated $X$
\end{tabular}
\caption{Comparison of augmented minimax linear estimation \eqref{eq:debias} and the estimator of
\citet{wooldridge2018inference} \eqref{eq:debiasWZ} for average partial effect estimation with a logistic link. The
boxplots depict \smash{$\htau_1$}-estimates across 20 simulation replications for each method; the dashed
line is the true average partial effect.}
\label{fig:simu}
\end{center}
\end{figure}

As seen in Figure \ref{fig:simu}, both methods perform reasonably well when $X$ is uncorrelated.  Augmented minimax
linear estimation is somewhat less variable in small samples; however, this may be due to the choice of tuning parameters (we used our
own implementation of the method of \citet{wooldridge2018inference}). When $X$ is correlated, both methods struggle; and this is a difficult
problem, as $\theta$ is not particularly sparse, and $X$ is correlated in a way that can induce confounding. Overall, however,
we see that the augmented minimax linear estimator is converging as $n$ increases, whereas the method of
\citet{wooldridge2018inference} is noticeably biased here. Thus, as reflected by the weaker assumptions required by our formal results,
augmented minimax linear estimation may be more robust to confounding in problems of this type.

\paragraph{Proof of Theorem \ref{theo:main}.}
We start by showing that the first term in \eqref{eq:error-expansion} is $o_p(n^{-1/2})$.
Our argument, which is a variant on one used in the proof of \citet[Theorem 2]{hirshberg2017balancing}, relies on the characterization
\begin{align*}
&I\p{\gamma; \htheta} = \sup_{f \in \ff_n} \frac{1}{n} \sum_{i = 1}^n  \sqb{h(X_i,f) - \gamma_i f(X_i) } \ \text{ for } \\ 
&h(x,f) = \frac{\partial f}{\partial x_j}(x),\ \ff_n = \{ \psi(x^T \htheta_n) x^T \theta : \norm{\theta}_{1} \le 1\}.
\end{align*}
With $\gamma_i = \gamma_i^{\star} = \gamma^{\star}_n(X_i)$, this average is centered for all functions $f \in \ff_n$, as a straightforward calculation shows that 
$\gamma^{\star}_n$ satisfies $\EE[]{ h(X_i,f) } = \EE[]{\gamma(X_i) f(X_i)}$ for all $f \in \ff_n$. Thus, 
$I(\gamma^{\star}; \htheta_n)$ is $O_p(R_n(\hh_n))$ where $R_n(\hh_n)$ is the Rademacher complexity of the class 
\[\begin{split}
&\hh_n = \set*{ h(x, f) - \gamma^{\star}_n(x)f(x) : f \in \ff_n } = \set*{ v(x)^T \theta :  \norm{\theta}_{1} \le 1} \ \text{ where }\\
& v(x) = \htheta_j \psi'\p{x^T \htheta} x + \psi\p{x^T \htheta} e_j -  \gamma^{\star}_n(x) \psi\p{x^T \htheta} x. 
\end{split}\]
To bound $R_n(\hh_n)$, we observe that $\hh_n$ is the convex hull of the finite set 
\[ \hh_n' = \set*{ \htheta_j \psi'\p{x^T \htheta} x_{k} + \psi\p{x^T \htheta} 1_{\set{k=j}} -  \gamma^{\star}_n(x) \psi\p{x^T \htheta} x_{k} : k \in 1 \ldots p_n }, \]
and thus $R_n(\hh_n) = R_n(\hh_n')$  \citep[see e.g.][Theorem 12]{bartlett2002rademacher}. We use Massart's Finite Class Lemma \citep[Lemma 5.2]{massart2000some}
to bound this quantity: $R_n(\hh_n') \le r_n \sqrt{2\log(p_n)}/n$ where $r_n^2 = \max_{h \in \hh_n'} E h(x)^2$. As $r_n = O(1)$ under our boundedness assumptions and
$\log(p_n) = o(n)$, it follows that $R_n(\hh_n') = o(n^{-1/2})$  and therefore that $I(\gamma^{\star}; \htheta_n) = o_p(n^{-1/2})$. 
Then as simple consequence of the criterion \eqref{eq:debias} we use to choose our weights and our assumption $\norm{\gamma_n^{\star}(\cdot)}_{\infty} = O(1)$, 
\[ I(\hgamma; \htheta)^2 \le I(\gamma^{\star}; \htheta)^2 + n^{-2}(\norm{\gamma^{\star}}^2 - \norm{\hgamma}^2) \le o_p(n^{-1}) + O(n^{-1}), \]
so we have $I(\hgamma; \htheta) = O_p(n^{-1/2})$.\footnote{ In fact, as a consequence of the convergence of $\hgamma$ to $\gamma^{\star}$,
which we will establish below, this criterion will imply that $I(\hgamma; \htheta) = o_p(n^{-1/2})$.}\ 
As the first term in \eqref{eq:error-expansion} is $n^{-1}\sum_{i=1}^n v(X_i)^T (\htheta_n - \theta_n)$,
it is bounded by $I(\hgamma; \htheta) \norm{\htheta - \theta_n}_{1}$. Given our assumption that $\htheta$ is $\ell_1$-consistent, 
this implies that this term is $o_p(n^{-1/2})$.

Our second step will be to show that the weights $\hgamma$ solving \eqref{eq:balance} converge to the weights $\gamma^{\star}$
in empirical mean-square. This is sufficient to establish that the second term in \eqref{eq:error-expansion} is 
$-n^{-1}\sum_{i=1}^n \gamma_n^{\star}(X_i) (Y_i - \Psi(X_i^T \htheta_n)) + o_p(n^{-1/2})$. As the argument that this convergence property is sufficient appears
in the proof of \citet[Theorem 2]{hirshberg2017balancing}, we will not repeat it here. \citet[Theorem 4]{hirshberg2017balancing} establishes this convergence
under its conditions (i)-(vi), so it suffices to show that they are satisfied. In the language of that theorem, we will take $\tilde{\ff_n} = \ff_{L,n} = \ff_n$.
To save space, we will show that these conditions are satisfied without restating them here. 

Condition (i) is, stated more concretely, the continuity of the function $\theta \to (\htheta_j \psi'(X_i^T \htheta) X_i^T + \psi(X_i^T \htheta) e_j^T)\theta$ as a map from $\ell_1$ to $\ell_1$.
This is implied by our those boundedness assumptions on $X_i, \psi, \psi'$. 
Conditions (ii) and (iv) follow from our boundedness assumptions on these quantities and on $\gamma^{\star}_n$.
Condition (iii), in which we take $\tilde{\gamma}_n$ to be $\gamma^{\star}_n$, follows from our assumption that $g_n(\htheta_n) = o(n^{1/2})$. 
And while the parts of condition (vi) involving $\hm - m$ are not satisfied, 
the only part that is used in the proof of \citet[Theorem 2]{hirshberg2017balancing} to show convergence of the weights is the 
condition $R_n(\hh_n) = O_p(n^{-1/2})$ that we established above.\footnote{The parts of condition (vi) involving $\hm - m$ are used
only to control the first term in our error expansion, which we treated above using a variation on the argument used in \citet{hirshberg2017balancing}.}\ 
This leaves condition (v), which is stated in terms of 
\[ \ff_n^{\star}(r) = (\ff_n - [0,1]\gamma^{\star}_n)\cap rB \ \text{ and } \hh_n^{\star}(r) = \set*{ h(x, f) - \gamma^{\star}_n(x)f(x) : f \in \ff_n(r) }  \]
where $B$ is the unit $L_2$ ball $\set{ f : E f(X)^2 \le 1 }$. As $\ff_n^{\star}(r) \subseteq \ff_n^{\star}(\infty)$ and $\hh_n^{\star}(r) \subseteq \hh_n^{\star}(\infty)$,
condition (v) is implied by the bounds $R_n(\ff_n^{\star}(\infty)) = o(n^{-1/2})$ and $R_n(\hh_n^{\star}(\infty)) = o(n^{-1/2})$. 
We will show only the first of these bounds, as the argument for the second is analogous. 
Because $R_n(\ff + \ff') \le R_n(\ff) + R_n(\ff')$ for any sets $\ff,\ff'$ \citep[see e.g.][Theorem 12]{bartlett2002rademacher},
$R_n(\ff_n^{\star}(\infty)) \le R_n(\ff_n) + R_n(-[0,1]\gamma^{\star}_n) = R_n(\ff_n') + R_n(\set{0,-\gamma^{\star}_n})$ where
$\ff_n' = \set{ \psi(x^T \htheta_n) x_k : k \in 1 \ldots p_n}$ and $\set{0,-\gamma^{\star}_n}$ are finite classes which have $\ff_n$ and $-[0,1]\gamma^{\star}_n$ as their respective convex hulls. As the elements of these finite classes are bounded uniformly in $n$, $R_n(\ff_n') = O(\sqrt{\log(p_n)}/n) = o(n^{-1/2})$ and $R_n(\set{0,-\gamma^{\star}_n}) = O(n^{-1})$
by Massart's finite class lemma. Thus, we have established our claimed bound and therefore (v). This completes our proof of our asymptotic characterization of
the second term in \eqref{eq:error-expansion}.

Our final step will be to show that the remainder $rem_n$ in \eqref{eq:error-expansion} is $o_p(n^{-1/2})$. Subtracting the first expression in \eqref{eq:error-expansion} from the 
second shows that $-rem_n = n^{-1}\sum_{i=1}^n a(X_i, \hgamma_i)$ where
\begin{align*}
a(x,\gamma) &= \theta_j \psi\p{x^T \theta} - \htheta_j \psi\p{x^T \htheta} -	
	\p{\htheta_j \psi'\p{x^T \htheta} x^T + \psi\p{x^T \htheta} e_j^T} \p{\theta - \htheta} \\
& - \gamma \p{\Psi\p{x^T \theta} - \Psi\p{x^T \htheta} - \psi\p{x^T \htheta} x^T\p{\theta - \htheta}}.
\end{align*}
By design, $a(x,\gamma)=b(x,\gamma, \theta)-b(x,\gamma, \htheta)-(\nabla_{\theta}\mid_{\theta=\htheta} b(x,\gamma,\theta)) (\theta-\htheta)$ for
$b(x,\gamma, \theta) = \theta_j \psi(x^T \theta) -  \gamma \Psi(x^T \theta)$. It follows that
$-rem_n = b(\theta) - b(\htheta) - (\nabla_{\theta}\mid_{\theta=\htheta} b(\theta)) (\theta-\htheta)$ for $b(\theta)=n^{-1}\sum_{i=1}^n b(X_i,\hgamma_i,\theta)$,
i.e. $-rem_n$ is the remainder after first-order Taylor approximation of this function $b$ around $\htheta$ evaluated at $\theta$. 
Thus, using the Lagrange form of the remainder after Taylor approximation, we have $-rem_n= (1/2) (\theta - \htheta)^T H(\ttheta) (\theta - \htheta)$
where $H(\ttheta)$ is the Hessian of $b$ at some vector $\ttheta$ on the line segment between $\theta$ and $\htheta$,
\[ H(\ttheta) = n^{-1}\sum_{i=1}^n \sqb{\psi'\p{X_i^T \ttheta}(X_i e_j^T + e_j^T X_i) + \ttheta_j \psi''\p{X_i^T \ttheta} X_i X_i^T - \hgamma_i \psi'\p{X_i^T \ttheta} X_i X_i^T}, \]
Letting $z_n=\theta_n-\htheta_n$, we write $z_n^T H(\ttheta) z_n$ as a sum of three terms $\xi_{1,n}(\ttheta) + \xi_{2,n}(\ttheta) + \xi_{3,n}(\ttheta)$ for
\begin{align*}
&\xi_{1,n}(\ttheta) = 2 z_{j,n} n^{-1}\sum_{i=1}^n \psi'\p{X_{i,n}^T \ttheta_n} X_{i,n}^T z_n; \\
&\xi_{2,n}(\ttheta) = \ttheta_{j,n} n^{-1}\sum_{i=1}^n \psi''\p{X_{i,n}^T \ttheta_n} (X_{i,n}^T z_n)^2; \\
&\xi_{3,n}(\ttheta) = -n^{-1}\sum_{i=1}^n \hgamma_i \psi'\p{X_{i,n}^T \ttheta_n} (X_{i,n}^T z_n)^2. 
\end{align*}
We can bound $\abs{\xi_{k,n}(\ttheta)}$ by $\sup_{t \in [0,1]}\abs{\xi_{k,n}(\ttheta(t))}$ for $\ttheta(t)=\htheta + t (\theta - \htheta)$,
and by Markov's inequality this quantity will be $O_p(E \sup_{t \in [0,1]}\abs{\xi_{k,n}(\ttheta(t))})$, where
\begin{align*} 
E \sup_{t \in [0,1]} \abs*{\xi_{1,n}}
&=  E \sup_{t \in [0,1]}\abs*{ \psi'\p{X_{i,n}^T \ttheta_n} } z_n^T (e_j X_{i,n}^T + X_{i,n}^T e_j^T) z_n \\
&\le \norm{ \psi' }_{\infty} z_n^T A_1 z_n \ \text{ for } A_1 = E e_j X_{i,n}^T + X_{i,n}^T e_j^T ; \\
E \sup_{t \in [0,1]} \abs*{\xi_{2,n}}
&=  E \sup_{t \in [0,1]}\abs*{ \ttheta_{j,n} \psi''\p{X_{i,n}^T \ttheta_n} } z_n^T X_{i,n} X_{i,n}^T z_n \\
&\le \max\set{\abs{\theta_{j,n}},\abs{\htheta_{j,n}}} \norm{ \psi'' }_{\infty} z_n^T A_2 z_n \ \text{ for } A_2 = E X_{i,n} X_{i,n}^T.
\end{align*}
As we've assumed that $\psi'$ and $\psi''$ are bounded and that $\theta_{j,n}=O(1)$ and this latter assumption and our assumption that $\norm{\htheta_n - \theta_n}_1 =o(1)$
implies that $\htheta_{j,n}=O(1)$, it follows that these quantities are $O(z_n^T A_k z_n)$.
We've assumed that $z_n^T A_k z_n = o(n^{-1/2})$ for $k \in \set{1,2}$,
so it follows that $\xi_{k,n}(\ttheta) = o_p(n^{-1/2})$ for $k \in \set{0,1}$.

Because $\hgamma_i$ is dependent on $X_1 \ldots X_n$, we cannot use this argument directly to bound $\xi_{3,n}(\ttheta)$. 
To work around this, we will first show that both $\xi_{3,n}(\ttheta)-\xi_{3,n}'(\ttheta)$ and $\xi_{3,n}'(\ttheta)$ are $o_p(n^{-1/2})$ where
\[ \xi_{3,n}'(\ttheta) = -n^{-1}\sum_{i=1}^n \gamma_n^{\star}(X_{i,n}) \psi'\p{X_{i,n}^T \ttheta} (X_{i,n}^T z)^2. \] 
To bound $\xi_{3,n}(\ttheta)-\xi_{3,n}'(\ttheta)$, we use the Cauchy-Schwartz inequality,
\begin{align*}
\xi_{3,n}'(\ttheta)-\xi_{3,n}(\ttheta) &= n^{-1}\sum_{i=1}^n \p{\hgamma_i - \gamma_n^{\star}(X_{i,n})} \psi'\p{X_{i,n}^T \ttheta} (X_{i,n}^T z_n)^2 \\
					 &\le \sqrt{ n^{-1}\sum_{i=1}^n \p{\hgamma_i - \gamma_n^{\star}(X_{i,n})}^2} \cdot \sqrt{n^{-1}\sum_{i=1}^n \p{\psi'\p{X_{i,n}^T \ttheta}^2 (X_{i,n}^T z_n)^2}^2 }.
\end{align*}
In this bound, the first factor is $o_p(1)$ as a consequence of the mean-square consistency of $\hgamma$. The second factor is amenable to the approach we've used above to bound
$\xi_{k,n}(\ttheta)$ for $k \in \set{1,2}$, which shows that $\xi_{3,n}(\ttheta)-\xi_{3,n}'(\ttheta) = O_p(z_n^T A_2 z_n) = o_p(n^{-1/2})$. 
And as we've assumed that $\norm{\gamma_n^{\star}}_{\infty} = O(1)$, the same argument yields a bound 
$\xi_{3,n}'(\ttheta) = O_p(z_n^T A_2 z_n) = o_p(n^{-1/2})$. This completes our proof that $rem_n$ is $o_p(n^{-1/2})$.

Using our characterizations of all three terms in \eqref{eq:error-expansion}, we have
$\htau_j - \htau^{\star}_j = n^{-1}\sum_{i=1}^n \gamma_n^{\star}(X_i) (Y_i - \Psi(X_i^T \htheta_n)) + o_p(n^{-1/2})$.
Adding $\htau^{\star}_j - \tau_j$ yields our claimed asymptotic characterization \eqref{eq:asymptotic-linear}.

\bibliographystyle{plainnat}
\bibliography{references}

\end{document}